\documentclass{amsart}
\newtheorem{teor}{Theorem}[section]
\newtheorem{prop}[teor]{Proposition}
\newtheorem{coro}[teor]{Corollary}
\newtheorem{lema}[teor]{Lemma}
\newtheorem{defi}[teor]{Definition}

\newtheorem{nota}{Remark}[section]

\def\qed{\ \ \ \hbox{}\nolinebreak\hfill $\blacksquare
 \  \  \  \  $ \par{}\medskip}

\def\ba{\begin{eqnarray*}}
\def\ea{\end{eqnarray*}}

\def\be{\begin{equation}}
\def\ee{\end{equation}}

\def\bd{\begin{defi}}
\def\ed{\end{defi}}

\def\bp{\begin{prop}}
\def\ep{\end{prop}}

\def\bt{\begin{teor}}
\def\et{\end{teor}}

\def\bl{\begin{lema}}
\def\el{\end{lema}}

\def\bn{\begin{nota}}
\def\en{\end{nota}}

\def\bc{\begin{coro}}
\def\ec{\end{coro}}

\newfont{\sssB}{msam10 scaled\magstephalf}
\def\sqb{{\sssB\char4}}

\def\ds{\displaystyle}

\def\m{\vert}

\def\qed{\sqb \medskip}

\def\md1{\sum_{k=1}^n \m x_k\m}

\newfont{\bbb}{msbm10 scaled\magstephalf}

\def\R{\mbox{\bbb R}}
\def\Z{\mbox{\bbb Z}}
\def\N{\mbox{\bbb N}}

\def\T{\mbox{\bbb T}}

\def\Dem{{\it Proof. }}

\def\ds{\displaystyle}

\def\lpq{{L^{p,q}}}

\def\lppqq{{L^{p',q'}}}

\def \ecu{\hskip 5pt }

\def \ep{\epsilon}

\def \al{\alpha }

\def \toro{{\rm T\! \! \! T}}

\def \num{${\mbox n}^{\underline{\mbox{\tiny o}}}$ } 
\def \12{\frac{1}{2}}

\def \dint {\int \! \! \! \int }

\def\bbbr{{\rm I\!R}} %reelle Zahlen

\def\bbbz{{\rm Z\!\!Z}}
\def\bbbn{{\rm I\!N}} %natuerliche Zahlen

\def \num{${\mbox n}^{\underline{{\rm \tiny o}}}$}

\begin{document}

%%%%%%%%%%%%%%%%%%%%%%%%%%%%%%%%%%%

\centerline{AMS Class (2000) 42A45}

\begin{abstract} Let $m(\xi,\eta)$ be a bounded continuous function in $\bbbr\times\bbbr$, $0<
p_i,q_i<\infty$ for $i=1,2$ and $0<p_3,q_3\le\infty$ 
  where
$1/p_1+1/p_2=1/p_3$. It is
shown that 
$$C_m (f,g)(x)=\int_{\bbbr} \int_{\bbbr} \hat f(\xi) \hat g(\eta) m(\xi,\eta) e^{2\pi i x(\xi +\eta
)}d\xi
d\eta$$ is a bounded bilinear operator from $L^{p_1,q_1}(\bbbr)\times
L^{p_2,q_2}(\bbbr)$ into
$L^{p_3,q_3}(\bbbr)$  if and only
if 
$$P_{D_{\varepsilon^{-1}}m} (f,g)(\theta)=\sum_{k\in \Z} \sum_{k'\in \Z}  \hat f(k) \hat g(k')
m(\varepsilon k,
\varepsilon k') e^{2\pi i
\theta(k +k' )}$$ are  bounded bilinear operators  from $L^{p_1,q_1}(\T)\times
L^{p_2,q_2}(\T)$ into
$L^{p_3,q_3}(\T)$ with norm bounded by uniform constant for all $\epsilon >0$.
\end{abstract}

\title{Transference of bilinear multiplier operators on Lorentz spaces.}

\author{Oscar Blasco, Francisco Villarroya}
\thanks{Both authors have been partially supported by grants DGESIC 
PB98-0146 and BMF 2002-04013}
\date{}
\maketitle

%%%%%%%%%%%%%%%%%%%%%%%%%%%%%%%%%%%%%%%%%%%%%%%%%%%%%%%%%
%%%%%%%%%%%%%%%%%%%%%%%%%%%%%%%%%%%%%%%%%%%%%%%%%%%%%%%%%

\section{Introduction.}

Let $m(\xi_1,\xi_2,...,\xi_n)$ be a bounded measurable function in $\bbbr^n$ and define
 $$C_m (f_1,f_2,...,f_n)(x)=\int_{\bbbr^n} \hat f_1(\xi_1)... \hat f_n(\xi_n) m(\xi_1,\xi_2,...,\xi_n)
e^{2\pi i x(\xi_1 +\xi_2+...+\xi_n )}d\xi$$  for 
Schwartz test functions $f_i$ in ${\mathcal S }$ for $i=1,...,n$.

Given now $0<p_i\le\infty$ for $i=1,...,n$ and $1/q=1/p_1+1/p_2+...1/p_n$. The function $m$ is said to
be a multilinear
multiplier of strong type
$(p_1,p_2,...,p_n)$ (respect. weak type $(p_1,p_2,...,p_n)$) if $C_m$ extends to a bounded bilinear
operator from
$L^{p_1}(\bbbr)\times... \times L^{p_n}(\bbbr)$ into
$L^{q}(\bbbr)$ (respect. to $L^{q,\infty}(\bbbr)$).

The study of such multilinear multipliers   was started by R. Coifman and Y. Meyer (see \cite  {CM1,
CM2, CM3})
for smooth symbols. 
However, in  the last years people got interested in them  after the    results proved by  M.
Lacey and C. Thiele 
(\cite{LT1, LT2, LT3}) which establish that
$m(\xi,\nu)=sign(\xi+\alpha\nu)$ are multipliers of
strong type
$(p_1,p_2)$ for $1<p_1,p_2\le\infty$,
$p_3>2/3$ and each $\alpha\in \bbbr\setminus\{ 0,1\}$.

New results for non-smooth symbols, extending the ones given by the bilinear Hilbert transform,  have
been achieved  by J.E.
Gilbert and A.R. Nahmod (see \cite {GN1, GN2, GN3}) and  by C. Muscalu, T. Tao and C.
Thiele (see
\cite{MTT}).

We refer the reader to \cite{KG, KS, FS, GT} for several results on bilinear multipliers
and related topics.

The first transference methods for linear multipliers were given by K. Deleeuw.
It is known that if $m$ is continuous then  $$T_m (f)(x)= \int_{\bbbr} \hat f(\xi)  m(\xi) e^{2\pi i
x\xi 
}d\xi$$ 
(defined for $f\in S(\bbbr)$) is bounded on $L^p(\bbbr)$ if and
only if
$$\tilde T_{m_\varepsilon} (f)(\theta)=\sum_{k\in \Z}   \hat f(k)  m(\varepsilon k) e^{2\pi i \theta k
}$$ (defined for trigonometric polynomials $f$) are uniformly bounded on $L^p(\T)$ for all
$\varepsilon>0$ (see
\cite{L}, \cite{SW} page
264).

Although the results in the paper hold true for multilinear multipliers, for simplicity of the
notation we restrict ourselves
to bilinear multipliers and only state and prove the theorems in such a situation.

Let $(m_{k,k'})$ be a bounded sequence  we use the notation
$$P_{m} (f,g)(\theta)=\sum_{k\in \bbbz} \sum_{k'\in \bbbz}  a(k) b(k') m_{ k, k'} e^{2\pi i \theta(k
+k' )}$$
for 
$f(t)=\sum_{n\in\bbbz}
a(n)e^{2\pi i nt}$ and
$g(t)=\sum_{ n\in\bbbz} b(n)e^{2\pi i nt}$.

Let $0<p_1,p_2\le \infty$ and $p_3$ such that $1/p_1+1/p_2=1/p_3$. We write $P_{D_{t^{-1}}m}$ 
when the
symbol is $m(tk,tk')$ and
say that
$m(tk,tk')$ is a bounded multiplier of strong (respect. weak) type
$(p_1,p_2)$ on
$\bbbz\times\bbbz$ if the corresponding $P_{D_{t^{-1}}m}$ is bounded from
$L^{p_1}(\T)\times L^{p_1}(\T)$ into
$L^{p_3}(\T)$ (respect.
$L^{p_3,\infty}(\T)$).

 In
a recent paper (see \cite{FS})  D. Fan and S. Sato  have shown certain DeLeeuw type theorems for
transferring
multilinear operators on Lebesgue and Hardy spaces from $\bbbr^n$ to $\T^n$. 
They show that the multilinear version of the transference between $\bbbr$ and
$\bbbz$ holds true, namely
that for continuous functions $m(\xi,\eta)$ one has that $m$ is a multiplier of strong (respect. weak)
type $(p_1,p_2)$  on
$\bbbr\times\bbbr$ if and only if
$(D_{\varepsilon^{-1}}m)_{k,k'}=(m(\varepsilon k,\varepsilon k'))_{k,k'}$ are uniformly
bounded multipliers of strong (respect.
weak) type
$(p_1,p_2)$ on
$\bbbz\times\bbbz$.

The first author (see \cite{Bl})  has shown a Deleeuw type theorem to transfer
   bilinear multipliers from $L^p(\bbbr)$ to  bilinear multipliers acting on $\ell_p(\bbbz)$.

The aim of this paper is to get an extension of those results in  \cite{FS} for bilinear multipliers
acting on Lorentz spaces
(see  \cite{FS}, Remark 3). 

We shall show that if  $m$ is a bounded continuous   function on  $\bbbr^{2}$
then $C_m$ defines a bounded  bilinear map from $L^{p_{1},q_{1}}(\bbbr )\times
L^{p_{2},q_{2}}(\bbbr )$ into $L^{p_{3},q_{3}}(\bbbr )$ 
if and only if the $P_{D_{t^{-1}}m}$, the restriction to $m(tk,tk')$ for
$k,k'\in\bbbz$,  define 
 bilinear maps from $L^{p_{1},q_{1}}(\T )\times
L^{p_{2},q_{2}}(\T )$ into $L^{p_{3},q_{3}}(\T )$  uniformly bounded for $t>0$.

Throughout the paper $|A|$  denotes the Lebesgue measure of $A$ and we identify functions $f$
on $\T $ and periodic functions on $\bbbr $  with period $1$ defined on $[-\12 ,\12 )$, that is
$f(x)=f(e^{2\pi ix})$  and
$
\int_{\toro }f(z)dm (z)=\int_{-\12 }^{\12 }f(t)dt.
$
For  $0< p\le\infty$, we write
$
D_{t}^{p}f(x)=t^{-\frac{1}{p}}f(t^{-1}x)$ (with the notation $D_t=D^\infty_t$), $M_{y}f(x)=f(x)e^{2\pi
iyx}
$ and
$
T_{y}f(x)=f(x-y)
$ for the dilation, modulation and translation operators. In this way
$(D_{t}^{q}f)\hat{}=D_{t^{-1}}^{q'}\hat{f}
$ where, as usual, $q'$ stands for the conjugate exponent of $q$. 

Adknowledgement: We want to thank the referee for his or her carefull reading. 

\section{Preliminaries}

Let $(\Omega,
\Sigma,\mu)$ be a $\sigma$-finite and complete measure space.
 Given a complex-valued measurable function $f$ we shall  denote  the distribution function
of
$f$ by $\mu_f(\lambda) = \mu(E_\lambda)$ for $\lambda > 0$ where $E_\lambda = \{ w
\in
\Omega : |f(w)| > \lambda\}$, 
 the nonincreasing rearrangement of
$f$ by $f^*(t) = \inf\{\lambda>0 : \mu_f(\lambda) \leq t\}$ and  $f^{**}(t) =
\frac{1}{t} \int_0^t f^*(s) ds$.

Now the Lorentz space $\lpq$ consists of those
 measurable functions $f$ such that $\|f\|^*_{pq} < \infty$, where
$$
\|f\|^*_{pq} = \left\{
\begin{array}{ll}
\ds{\left\{ \frac{q}{p} \int_0^\infty t^{\frac{q}{p}} f^*(t)^q \frac{dt}{t}\right\}^\frac{1}{q},} &
0<p<\infty ,\ 0 < q < \infty , \\
\ds{\sup_{t> 0} t^{\frac{1}{p}} f^*(t) }  & 0 < p \leq \infty , \ q=\infty .
\end{array} \right.
$$

It is well known that 
$$||f||_{p\infty}= \sup_{       \lambda>0} \lambda \mu_f(\lambda)^{1/p}. $$
Here we shall use the following fact:
If $0<p,q<\infty$ and $f$ is a measurable function then
\begin{equation}\label{equiv}
\|f\|^*_{pq}=\left(
q \int_0^\infty \lambda^{q-1} \mu_f(\lambda)^{\frac{q}{p}} d\lambda\right)^{1/q}.
\end{equation}
(This can be  easily checked for simple functions).

Let us recall some facts about these spaces. Simple functions are dense in 
$\lpq$ for $q \neq \infty$,  
$(L^{p,1})^* = L^{p',\infty}$ for $1 \leq p < \infty,$ and
$(\lpq)^* =
\lppqq$ for $1 < p,q <
\infty$ as well.
 Replacing $f^*$ by $f^{**}$ and putting
$\|f\|_{pq} =
\|f^{**}\|^*_{pq}$ then we get a functional
 equivalent  to $\| \cdot \|^*_{pq}$ (for $1 < p < \infty$) for which  $L^{1,1}$ and $\lpq$ for $1
< p \leq \infty$, $1 \leq q \leq \infty$ are Banach spaces.

The reader is referred to \cite{H}, \cite{BS}, \cite{SW} or \cite {Lo} for
the basic information on Lorentz spaces. We only condider $\mu$ to be either the Lebesgue measure on
$\bbbr$ or the normalized
Lebesge measure on $\T$ and the distribution function will be denoted $m_f$ in both cases.

\begin{defi}
Let $m$ be a bounded measurable function on $\bbbr^{2}$. Let $0<p_{i},q_{i}\le\infty$ for $i=1,2,3$.
 For $t>0$ we define
$$
C_{D_{t^{-1}}m}(f,g)(x)=C_t(f,g)(x)
=\int_{\bbbr }\int_{\bbbr }\hat{f}(\xi )\hat{g}(\eta )m(t\xi ,t\eta )
e^{2\pi i(\xi +\eta )x}d\xi d\eta 
$$
for $f,g\in S(\bbbr )$ .

 We say that $m$ is a bilinear multipier in
$(L^{p_1,q_1}(\bbbr)\times L^{p_2,q_2}(\bbbr),L^{p_3,q_3}(\bbbr))$ 
if there exists $C>0$ such that$$
\| C_{1}(f,g)\|_{L^{p_{3},q_{3}}(\bbbr )}\leq C\|f\|_{L^{p_{1},q_{1}}(\bbbr )}
\|g\|_{L^{p_{2},q_{2}}(\bbbr )}
$$
for all $f,g\in S(\bbbr )$ .
\end{defi}

\begin{defi}
Let $(m_{k_1,k_2})_{k_1\in\bbbz,k_2\in\bbbz}$ be a bounded sequence.  Let $p_{i},q_{i}>0$ such that
$p_{3}^{-1}=p_{1}^{-1}+p_{2}^{-1}$. 
We define
$$
P_{m}(f,g)(x)
=\sum_{k_{1}\in \bbbz }\sum_{k_{2}\in \bbbz }a_{k_{1}}b_{k_{2}}m_{k_1,k_2}
e^{2\pi i(k_{1}+k_{2})x}
$$
for all trigonometric polynomials 
$f(x)=\sum_{|k|\leq N}a_{k}e^{2\pi ikx}$, $g(x)=\sum_{|k|\leq M}b_{k}e^{2\pi ikx}$ and $N,M\in \bbbn$.

We say that $m_{k,k'}$ is a bilinear multiplier in
 $(L^{p_1,q_1}(\T)\times L^{p_2,q_2}(\T),L^{p_3,q_3}(\T))$ 
if there exists $C>0$ such that$$
\| P_{m}(f,g)\|_{L^{p_{3},q_{3}}(\toro)}\leq C\|f\|_{L^{p_{1},q_{1}}(\toro )}
\|g\|_{L^{p_{2},q_{2}}(\toro )}
$$
for all trigonometric polynomials $f$ and $g$.

\end{defi}

\begin{nota}  $m$ is a multiplier in
 $(L^{p_1,q_1}(\bbbr)\times L^{p_2,q_2}(\bbbr),L^{p_3,q_3}(\bbbr))$ if and only if
$D_{t^{-1}}m(\xi,\eta)=m(t\xi,t\eta)$ is
also a
 multiplier  for each $t>0$.

Note that for each $t>0$ we have $m_{D_tf}(\lambda)=tm_f(\lambda)$. 
Hence
\begin{equation}
||D_t f||_{L^{p,q}(\bbbr)}=t^{1/p}||f||_{L^{p,q}(\bbbr)}. 
\end{equation}
for $0< p,q\le \infty$.  

Now the remark 
 follows easily from the formula
$$
C_{t}(f,g)=D_{t}C_{1}(D_{t^{-1}}f,D_{t^{-1}}g).
$$
Actually we have
$
\| C_{t}\|= \| C_{1}\| 
$ for all $t>0$.

\end{nota}

Let us start by recalling some facts to be used in the sequel.

\begin{defi} 
If $f$  is a measurable function on $\bbbr$ such that
$\max\{|f(x)|,|\hat f(x)|\}\le A/(1+|x|)^\alpha$ for some $A>0$ and
$\alpha>1$ then
$\tilde f$ stands for the well-defined periodic function (see
\cite{SW}, pages 250-253)
$$
\tilde f(x)=\sum_{k\in \bbbz }f(x+k)=\sum_{k\in \bbbz }\hat{f}(k)e^{2\pi ikx}.
$$
\end{defi}

\begin{lema}\label{realtoro} Let $0< p<\infty$ and 
$0<q\le\infty$. If $f\in S(\bbbr )$ we have
$$
4^{-\frac{1}{r}}\| f\|_{L^{p,q}(\bbbr )}\leq
\liminf_{t\rightarrow 0}t^{-\frac{1}{p}}\| \widetilde{D_{t}f} \|_{ L^{p,q}(\toro )}
$$
$$
\limsup_{t\rightarrow 0}t^{-\frac{1}{p}}\| \widetilde{D_{t}f} \|_{ L^{p,q}(\toro )}
\leq 4^{\frac{1}{r}}\| f\|_{L^{p,q}(\bbbr )}
$$
where $r=\log_{2}^{-1}(2^{\frac{1}{p}+1}\max{(2^{\frac{1}{q}-1},1)})$
and $\widetilde{D_{t}f}(x)=\sum_{k\in \bbbz }D_{t}f(x+k)$ 
is defined on $\T$.
%is the exponent given by Aoki-Rolewicz theorem.
\end{lema}
\Dem  Assume first that $f$ has compact support.
For $t>0$ small enough we have
$supp(D_{t}f)\subset [-\12 ,\12 ]$. This gives that
$$
\widetilde{D_{t}f} \chi_{[-\12 ,\12 ]}=D_{t}f \chi_{[-\12 ,\12 ]}=D_{t}f 
$$
In particular, for such $t$ we have
$$
m_{\widetilde{D_{t}f}}(\lambda )=|\{ x\in [-\12 ,\12 ] / |D_{t}f(x)|>\lambda \}|
=| \{ x\in \bbbr / |f(t^{-1}x)|>\lambda \}|
=tm_{f}(\lambda)
$$
and
$$
(\widetilde{D_{t}f})^{*}(x)=D_{t}(f^{*})(x), \qquad x>0.
$$
Hence
$$
\| \widetilde{D_{t}f} \|_{L^{p,q}(\T )}^{q}
=\frac{q}{p}\int_{0}^{1}
\Big( x^{\frac{1}{p}}(\widetilde{D_{t}f})^{*}(x)\Big)^{q}\frac{dx}{x}
$$
$$
=\frac{q}{p}\int_{0}^{1}
\Big( x^{\frac{1}{p}}f^{*}(t^{-1}x)\Big)^{q}\frac{dx}{x}
=t^{\frac{q}{p}}\frac{q}{p}\int_{0}^{t^{-1}}
\Big( x^{\frac{1}{p}}f^{*}(x)\Big)^{q}\frac{dx}{x}
$$
and therefore 
$$
\lim_{t\rightarrow 0}t^{-\frac{1}{p}}\| \widetilde{D_{t}f} \|_{L^{p,q}(\T )}
=\lim_{t\rightarrow 0}\Big( \frac{q}{p}\int_{0}^{t^{-1}}
\Big( x^{\frac{1}{p}}f^{*}(x)\Big)^{q}\frac{dx}{x}\Big)^{\frac{1}{q}}
=\| f\|_{L^{p,q}(\bbbr )}.
$$

The case $q=\infty$ is simpler.

For the general case, take $f_n= f\chi_{[-n,n]}$.
Observe that, for $|x|<1/2$
\begin{eqnarray*}
\widetilde{D_tf}(x)-\widetilde{D_tf_n}(x)&=&\sum_{k\in\bbbz} f(t^{-1}(x+k))-f_n(t^{-1}(x+k))\\
&=&\sum_{|k+x|>tn} f(t^{-1}(x+k))\\
\end{eqnarray*}

Hence, for any $m>0$,  we have that 
$$|\widetilde{D_tf}(x)-\widetilde{D_tf_n}(x)|\le\sum_{|k+x|>tn} \frac{C_m}{(1+t^{-1}|x+k|)^m}\le
t^m\sum_{|k+x|>tn}
\frac{C_m}{|x+k|^m}\le C_mt^m.$$

This shows that, selecting   $m>1/p$, we have
$$\lim_{t\to 0}t^{-1/p}||\widetilde{D_tf_n}-\widetilde{D_tf}||_{L^{\infty}(\toro )}\le C_m \lim_{t\to
0} t^{m-1/p}=0.$$

Given $\ep >0$, choose $n\in\N$ such that $$(1-\ep)||f||_{L^{p,q}(\bbbr)}\le
||f_n||_{L^{p,q}(\bbbr)}\le
||f||_{L^{p,q}(\bbbr)}$$

Now since $\| \cdot \|_{L^{p,q}(\bbbr)}$ is a quasi-norm with constant 
$C=2^{\frac{1}{p}}\max (2^{\frac{1}{q}-1},1)$ 
we have by the Aoki-Rolewic theorem \cite{aoki}
that $\| \cdot \|_{L^{p,q}(\bbbr)}$ is equivalent to a $r$-norm, namely 
$| \cdot |$, for $r=\log_{2}^{-1}(2C)$. More precisely we have
$$
|f|\leq \| f\|_{L^{p,q}(\bbbr)}\leq 4^{\frac{1}{r}}|f| 
$$
and thus we can write a triangular inequality to the $r$-power in the 
following way
$$
\| f+g\|_{L^{p,q}(\bbbr)}^{r}\leq 4(\| f\|_{L^{p,q}(\bbbr)}^{r}+
\| g\|_{L^{p,q}(\bbbr)}^{r}).
$$
 Hence, using this triangular inequality for $||.||^r_{L^{p,q}(\toro)}$ for the
corresponding power
$r\le 1$, according to
different values of $p$ and $q$, and  the previous case we get the desired formula.
\qed

\begin{lema} \label{tororealdos}
Let $0<p,q\le\infty$, $\varphi =\chi_{[-\12 ,\12 ]}$, $f\in {L^{p,q}(\T )}$ and
 $k\in \bbbn $. Then
$$
\| f\|_{L^{p,q}(\toro )}
=\| fD_{k}^{p}\varphi \|_{L^{p,q}(\bbbr )}
$$
\end{lema}
\Dem 
Using that $f$ is periodic we get
\begin{eqnarray*}
m_{fD_{k}^{p}\varphi }(\lambda )
&=&| \{ x\in \bbbr : |f(x)k^{-\frac{1}{p}}\chi_{[-\12 ,\12 ]} (k^{-1}x)|>\lambda \}|
\\
&=&| \{ x\in [-\frac{k}{2},\frac{k}{2}] : |f(x)|>k^{\frac{1}{p}}\lambda \}|\\
&=&k| \{ x\in [-\frac{1}{2},\frac{1}{2}] : |f(x)|>k^{\frac{1}{p}}\lambda \}|
=km_{f}(k^{\frac{1}{p}}\lambda ).
\end{eqnarray*}
Hence
$$
(fD_{k}^{p}\varphi )^{*}(t)
=\inf \{ \lambda >0 : km_{f}(k^{\frac{1}{p}}\lambda )<t\}
$$
$$
=k^{-\frac{1}{p}}
\inf \{ \lambda >0 : m_{f}(\lambda )<k^{-1}t\}
=D_{k}^{p}f^{*}(t)=(D_{k}^{p}f)^{*}(t)
$$

Therefore
$$
\| fD_{k}^{p}\varphi \|_{L^{p,q}(\bbbr )}^{q}
=\frac{q}{p}\int_{0}^{\infty }t^{\frac{q}{p}}(fD_{k}^{p}\varphi )^{*}(t)^{q}
\frac{dt}{t}
$$
$$
=\frac{q}{p}\int_{0}^{\infty }t^{\frac{q}{p}}k^{-\frac{q}{p}}f^{*}(k^{-1}t)^{q}
\frac{dt}{t}
=\frac{q}{p}\int_{0}^{\infty }t^{\frac{q}{p}}f^{*}(t)^{q}\frac{dt}{t}
=\| f\|_{L^{p,q}(\T )}^{q}.
$$
\qed

\begin{lema}\label{tororeal1}
Let $0< p< \infty$ and $f\in L^{p,\infty}(\T )$. If $\varphi \in S(\bbbr )$ is radial and decreasing
then
$$
 \limsup_{\ep \rightarrow 0}
\| fD_{\ep^{-1}}^{p}\varphi \|_{L^{p,\infty}(\bbbr )}
\leq 
\| \varphi \|_{L^{p}(\bbbr )}\| f\|_{L^{p,\infty}(\T )}.
$$
\end{lema}
\Dem
Note that for each $\ep >0$ and $\lambda>0$ we have
\begin{eqnarray*}
&&|\{x\in \bbbr : |f(x)\varphi (\ep x)|>t\} |=|\{|x|\leq 2^{-1}\lambda\ep^{-1} :
|f(x)\varphi (\ep x)|>t\}|\\
&& +\sum_{n=0}^{\infty }|\{2^{n-1}\lambda\ep^{-1}<|x|\leq 2^{n}\lambda\ep^{-1} : |f(x)\varphi (\ep
x)|>t\}|\\
&\leq& |\{|x|\leq 2^{-1}\lambda\ep^{-1} : |f(x)|>t\varphi (0)^{-1}\}|\\
&&+\sum_{n=0}^{\infty }
|\{2^{n-1}\lambda\ep^{-1}<|x|\leq 2^{n}\lambda\ep^{-1} : |f(x)|>t\varphi (\lambda
2^{n-1})^{-1}\}|\\ &\leq& |\{|x|\leq 2^{-1}([\lambda\ep^{-1}]+1) : |f(x)|>t\varphi
(0)^{-1}\}|\\ &&+\sum_{n=0}^{\infty }
|\{2^{n-1}[\lambda\ep^{-1}]<|x|\leq 2^{n}([\lambda\ep^{-1}]+1) : |f(x)|>t\varphi
(\lambda 2^{n-1})^{-1}\}|\\
&=&([\lambda\ep^{-1}]+1)|\{x\in \T : |f(x)|>t\varphi (0)^{-1}\}|\\
&&+\sum_{n=0}^{\infty }
(2^{n+1}([\lambda\ep^{-1}]+1)-2^{n}[\lambda\ep^{-1}])
|\{x\in \T : |f(x)|>t\varphi (\lambda 2^{n-1})^{-1}\}|\\
&\leq& (\lambda\ep^{-1}+1)|\{x\in \T : |f(x)|>t\varphi (0)^{-1}\}|\\
&&+\sum_{n=0}^{\infty }
2^{n}(\lambda\ep^{-1}+2)|\{x\in \T : |f(x)|>t\varphi (\lambda 2^{n-1})^{-1}\}|.
\end{eqnarray*}

Hence we get

\begin{equation}
 m_{f D_{\ep^{-1}}\varphi}(t)
\leq (\lambda \ep^{-1}+1)m_f(t\varphi (0)^{-1})+
(\lambda \ep^{-1}+2)\sum_{n=0}^{\infty }
2^{n}m_f(t\varphi (\lambda 2^{n-1})^{-1}).
\end{equation}

Therefore, using that $m_f(t)\le \frac{||f||^p_{p\infty}}{t^p}$, we get 
\begin{eqnarray*}
m_{f D^p_{\ep^{-1}}\varphi}(s)&=&m_{f D_{\ep^{-1}}\varphi}(s\ep^{-1/p})\\
&\le&(\lambda\ep^{-1}+1)\ep s^{-p}\varphi (0)^{p}||f||^p_{L^{p,\infty}(\T)}\\
&+&\sum_{n=0}^{\infty }
2^{n}(\lambda\ep^{-1}+2)\ep s^{-p}\varphi (\lambda
2^{n-1})^{p}||f||^p_{L^{p,\infty}(\T)}\\ &\le&s^{-p}(\lambda+\ep) |\varphi
(0)|^{p}||f||^p_{L^{p,\infty}(\T)}\\ &+&s^{-p}\sum_{n=0}^{\infty }
2^{n}(\lambda+2\ep)\varphi (\lambda 2^{n-1})^{p}||f||^p_{L^{p,\infty}(\T)}
\end{eqnarray*}

Hence, if $\varphi_\lambda=\varphi (0)\chi_{[-\lambda 2^{-1},\lambda 2^{-1}]}
+\sum_{n\geq 0}\varphi (\lambda 2^{n-1})
\chi_{[-\lambda 2^{n},\lambda 2^{n}]\backslash [-\lambda 2^{n-1},\lambda 2^{n-1}]}$ we have
$$
\limsup_{\ep\to 0} ||f D^p_{\ep^{-1}}\varphi||^p_{L^{p,\infty}(\bbbr)}
\le
\|\varphi_\lambda \|^p_{L^{p}(\bbbr
)}||f||^p_{L^{p,\infty}(\T)}.
$$
Now pass to the limit as $\lambda$ goes to zero to get the result.
\qed

\begin{lema}\label{tororeal2}
Let $0< p,q< \infty$ and $f\in L^{p,q}(\T )$. If $\varphi \in S(\bbbr )$ is radial and decreasing then
$$
C_{p,s}
\| \varphi \|_{L^{p,s }(\bbbr )}\| f\|_{L^{p,q}(\T )}
\leq \liminf_{\ep \rightarrow 0}
\| fD_{\ep^{-1}}^{p}\varphi \|_{L^{p,q}(\bbbr )}$$
$$\leq \limsup_{\ep \rightarrow 0}
\| fD_{\ep^{-1}}^{p}\varphi \|_{L^{p,q}(\bbbr )}
\leq C_{p,r}
\| \varphi \|_{L^{p,r}(\bbbr )}\| f\|_{L^{p,q}(\T )}
$$
where $C_{p_1,p_2}=(2^{\frac{p_2}{p_1}}-1)^{-\frac{1}{p_2}}$, $r=\min (p,q)$ and $s=\max (p,q)$ . 
\end{lema}
\Dem Use (\ref{equiv}) to write 
\begin{eqnarray*}
\| fD_{\ep^{-1}}^{p}\varphi \|_{L^{p,q}(\bbbr )}^{q}
&=&\int_{0}^{\infty }qt^{q-1}( m_{fD_{\ep^{-1}}}(\ep^{-1/p}t))^{\frac{q}{p}}dt\\
&=&\int_{0}^{\infty }qt^{q-1}(\ep m_{fD_{\ep^{-1}}}(t))^{\frac{q}{p}}dt
\end{eqnarray*}
Using the estimate in the previous lemma we have
$$
\ep m_{f D_{\ep^{-1}}\varphi}(t)
\leq (\lambda + \ep  )m_f(t\varphi (0)^{-1})+
(\lambda+ 2\ep )\sum_{n=0}^{\infty }
2^{n}m_f(t\varphi (\lambda 2^{n-1})^{-1}).
$$

Now we see that for $r=\min (p,q)$ we have
\begin{equation}\label{esti}
\limsup_{\ep \rightarrow 0}\| fD_{\ep^{-1}}^{p}\varphi \|_{L^{p,q}(\bbbr )}
\leq\Big( \lambda^{\frac{r}{p}} \varphi (0)^{r}
+\sum_{n=0}^{\infty }(\lambda 2^{n})^{\frac{r}{p}}
\varphi (\lambda 2^{n-1})^{r}\Big)^{\frac{1}{r}} \| f\|_{L^{p,q}(\toro )}.
\end{equation}

If  $q\leq p$ then, for every $\lambda$, we have 
\begin{eqnarray*}
&&\| fD_{\ep^{-1}}^{p}\varphi\|_{L^{p,q}(\bbbr )}^{q}
=\int_{0}^{\infty }qt^{q-1}\Big( \ep |\{x\in \bbbr  
|f(x)\varphi (\ep x)|>t\}|\Big)^{\frac{q}{p}}dt
\\\
&\leq & \int_{0}^{\infty }qt^{q-1}
\Big((\lambda +\ep)m_f(t\varphi (0)^{-1})+
(\lambda +2\ep)\sum_{n=0}^{\infty }
2^{n}m_f(t\varphi (\lambda 2^{n-1})^{-1})
\Big)^{\frac{q}{p}}dt
\\
&\leq& \int_{0}^{\infty }qt^{q-1}
 (\lambda +\ep)^{\frac{q}{p}}m_f(t\varphi (0)^{-1})^{\frac{q}{p}}dt\\
&&+
\int_{0}^{\infty }qt^{q-1}(\lambda +2\ep)^{\frac{q}{p}}\sum_{n=0}^{\infty }
2^{n{\frac{q}{p}}}m_f(t\varphi (\lambda 2^{n-1})^{-1})^{\frac{q}{p}}dt
\\
&=&(\lambda+ \ep )^{\frac{q}{p}}\varphi (0)^{q}
\int_{0}^{\infty }qt^{q-1}
m_f(t)^{\frac{q}{p}}dt
\\
&&+(\lambda +2\ep )^{\frac{q}{p}}\sum_{n=0}^{\infty }2^{n\frac{q}{p}}
\varphi (\lambda 2^{n-1})^{q}\int_{0}^{\infty }qt^{q-1}
m_f(t)^{\frac{q}{p}}dt
\\
&=&\Big( (\lambda+ \ep )^{\frac{q}{p}}|\varphi (0)|^{q}
+(\lambda +2 \ep )^{\frac{q}{p}}\sum_{n=0}^{\infty }2^{n\frac{q}{p}}
\varphi (\lambda 2^{n-1})^{q}\Big)
\| f\|_{L^{p,q}(\toro )}^{q}
\end{eqnarray*}
Therefore 
$$
\limsup_{\ep \rightarrow 0}
\| fD_{\ep^{-1}}^{p}\varphi \|_{L^{p,q}(\bbbr )}
\leq \Big( \lambda^{\frac{q}{p}}\varphi (0)^{q}
+\sum_{n=0}^{\infty }(\lambda 2^{n})^{\frac{q}{p}}\varphi (\lambda
2^{n-1})^{q}\Big)^{\frac{1}{q}}
\| f\|_{L^{p,q}(\toro )},
$$
which gives  (\ref{esti}).

In the case $q>p$ we can use Minkowski and get
\begin{eqnarray*}
\| fD_{\ep^{-1}}^{p}\varphi \|_{L^{p,q}(\bbbr )}^{p}
&=&\Big( \int_{0}^{\infty }(q^{\frac{p}{q}}t^{p(1-\frac{1}{q})}\ep |\{x\in \bbbr : 
|f(x)\varphi (\ep x)|>t\}|)^{\frac{q}{p}}dt\Big)^{\frac{p}{q}}
\\
&\leq& \Big( \int_{0}^{\infty }
\Big( q^{\frac{p}{q}}t^{p(1-\frac{1}{q})}(\lambda +\ep)m_f(t\varphi
(0)^{-1})
\\
&+&(\lambda +2\ep)\sum_{n=0}^{\infty }2^{n}q^{\frac{p}{q}}t^{p(1-\frac{1}{q})}
m_f(t\varphi (\lambda 2^{n-1})^{-1})\Big)^{\frac{q}{p}}
dt\Big)^{\frac{p}{q}}
\\
&\leq& (\lambda +\ep)\Big( \int_{0}^{\infty }
\Big(q^{\frac{p}{q}}t^{p(1-\frac{1}{q})}m_f(t\varphi
(0)^{-1})\Big)^{\frac{q}{p}} dt\Big)^{\frac{p}{q}}
\\
&+&(\lambda +2 \ep )\sum_{n=0}^{\infty }
2^{n}\Big( \int_{0}^{\infty }\Big(q^{\frac{p}{q}}t^{p(1-\frac{1}{q})}
m_f(t|\varphi (\lambda 2^{n-1})|^{-1})\Big)^{\frac{q}{p}}
dt\Big)^{\frac{p}{q}}
\\
&=&(\lambda +\ep)\varphi (0)^{p}\Big( \int_{0}^{\infty }
qt^{q-1}m_f(t)^{\frac{q}{p}}
dt\Big)^{\frac{p}{q}}
\\
&+&(\lambda + 2\ep )\sum_{n=0}^{\infty }
2^{n}\varphi (\lambda 2^{n-1})^{p}\Big( \int_{0}^{\infty }qt^{q-1}
m_f(t)^{\frac{q}{p}}
dt\Big)^{\frac{p}{q}}
\\
&=&\Big( (\lambda +\ep)\varphi (0)^{p}
+(\lambda +2\ep)\sum_{n=0}^{\infty }2^{n}\varphi (\lambda 2^{n-1})^{p}\Big)
\| f\|_{L^{p,q}(\T )}^{p}
\end{eqnarray*}
Therefore 
$$
\limsup_{\ep \rightarrow 0}\| fD_{\ep^{-1}}^{p}\varphi \|_{L^{p,q}(\bbbr )}
\leq \Big( \lambda \varphi (0)^{p}
+\sum_{n=0}^{\infty }\lambda 2^{n}\varphi (\lambda 2^{n-1})^{p}\Big)^{\frac{1}{p}}
\| f\|_{L^{p,q}(\T )},
$$
and  (\ref{esti}) is proved.

If $\varphi_\lambda=\varphi (0)\chi_{[-\lambda 2^{-1},\lambda 2^{-1}]}+\sum_{n\geq 0}\varphi (\lambda
2^{n-1})
\chi_{[-\lambda 2^{n},\lambda 2^{n}]\backslash [-\lambda 2^{n-1},\lambda 2^{n-1}]}$ then
clearly we have that
$$||\varphi_\lambda||_p=\Big( \lambda \varphi (0)^{p}
+\sum_{n=0}^{\infty }\lambda 2^{n}\varphi (\lambda 2^{n-1})^{p}\Big)^{\frac{1}{p}}.$$

Since $\varphi $ and $\varphi_\lambda$ are radial and decreasing  then
$\varphi_\lambda^*(t)=\varphi_\lambda(2t)$ for $t>0$
and 
$$|| \varphi_\lambda||^*_{pr}=\Big(\lambda^{\frac{r}{p}} \varphi (0)^{r}
+(2^{\frac{r}{p}}-1)\sum_{n=0}^{\infty }(\lambda2^{n})^{\frac{r}{p}}
\varphi (\lambda 2^{n-1})^{r}\Big)^{\frac{1}{r}}.$$

Hence, using that $r\le p$, we have
$$
\Big(\lambda^{\frac{r}{p}} \varphi (0)^{r}
+\sum_{n=0}^{\infty }(\lambda2^{n})^{\frac{r}{p}}
\varphi (\lambda 2^{n-1})^{r}\Big)^{\frac{1}{r}}
\le(2^{\frac{r}{p}}-1)^{-\frac{1}{r}}
\|\varphi_\lambda \|_{L^{p,r}(\bbbr )}.
$$

Finally taking limits as $\lambda\to 0$ give
$$
\limsup_{\ep \rightarrow 0}\| fD_{\ep^{-1}}^{p}\varphi \|_{L^{p,q}(\bbbr )}
\leq \lim_{\lambda \rightarrow 0}\Big( \lambda^{\frac{r}{p}} \varphi (0)^{r}
+\sum_{n=0}^{\infty }(\lambda 2^{n})^{\frac{r}{p}}
\varphi (\lambda 2^{n-1})^{r}\Big)^{\frac{1}{r}}
$$
$$
\leq (2^{\frac{r}{p}}-1)^{-\frac{1}{r}}
\limsup_{\lambda \rightarrow 0}\| \varphi_\lambda\|_{L^{p,r}(\bbbr
)}=(2^{\frac{r}{p}}-1)^{-\frac{1}{r}}
\| \varphi \|_{L^{p,r}(\bbbr )}.
$$

This gives one of the inequalities of the Lemma.

To get the other inequality, 
we use  estimates from below to obtain
$$
\liminf_{\ep \rightarrow 0}\| fD_{\ep^{-1}}^{p}\varphi \|_{L^{p,q}(\bbbr )}
\geq \Big( \lambda^{\frac{s}{p}}\varphi (\lambda 2^{-1})^{s}
+\sum_{n=0}^{\infty }(\lambda 2^{n})^{\frac{s}{p}}\varphi (\lambda
2^{n})^{s}\Big)^{\frac{1}{s}}
\| f\|_{L^{p,q}(\toro )}
$$
where $s=\max (p,q)$.

 Using now that $s\ge p$, we get, arguing as above, that
$$
\Big( \lambda^{\frac{s}{p}}\varphi (\lambda 2^{-1})^{s}
+\sum_{n=0}^{\infty }(\lambda 2^{n})^{\frac{s}{p}}\varphi (\lambda
2^{n})^{s}\Big)^{\frac{1}{s}}
\geq (2^{\frac{s}{p}}-1)^{-\frac{1}{s}}
\|\varphi^\lambda\|_{L^{p,s}(\bbbr )}
$$
where $\varphi^\lambda= \varphi (\lambda  2^{-1})\chi_{[-\lambda  2^{-1},\lambda 2^{-1}]}+\sum_{n\geq
0}\varphi (\lambda2^{n})
\chi_{[-\lambda 2^{n},\lambda 2^{n}]\backslash [-\lambda 2^{n-1},\lambda 2^{n-1}]}$.

Hence
$$
\liminf_{\ep \rightarrow 0}\| fD_{\ep^{-1}}^{p}\varphi \|_{L^{p,q}(\bbbr )}
\geq  (2^{\frac{s}{p}}-1)^{-\frac{1}{s}}
\| \varphi \|_{L^{p,s}(\bbbr )}.
$$

Then proof is then completed.\qed

\begin{coro}\label{tororeal3}
Let $0<p< \infty$ and $f\in L^{p}(\T )$. If $\varphi \in S(\bbbr )$ is radial and decreasing then
$$
\| \varphi \|_{L^{p }(\bbbr )}\| f\|_{L^{p}(\toro )}
= \lim_{\ep \rightarrow 0}
\| fD_{\ep^{-1}}^{p}\varphi \|_{L^{p}(\bbbr )}
.$$

In particular for $p=1$ and the periodized function $f=\widetilde{\chi_{A}}$ where
$A\subset [-\12 ,\12 ]$ we get
$$
\lim_{\ep \rightarrow 0}\int_{\bbbr }f(x)D_{\ep^{-1}}^{1}\varphi (x)dx=
\lim_{\ep \rightarrow 0}\int_{\bbbr }D_{\ep }f(x)\varphi (x)dx
=m (A)\int_{\bbbr}\varphi (x)dx
.$$
\end{coro}

%\begin{lema}Let $f\in L^{\infty }(\toro )$ and $\varphi \in L^{\infty }(\bbbr )$
%continuous at the origin and such that $\| \varphi \|_{\infty }=\varphi (0)=1$. Then 
%$$
%\| f\|_{L^{\infty }(\bbbr )}=\lim_{\ep \rightarrow 0}
%\| fD_{\ep^{-1}}^{1}\varphi \|_{L^{\infty }(\bbbr )}
%$$
%\end{lema}
%\Dem For a.e. $x\in \bbbr $ we have 
%$$
%|f(x)D_{\ep^{-1}}^{1}\varphi (x)|=|f(x)\varphi (\ep x)|
%\leq \| f\|_{L^{\infty }(\bbbr )}\| \varphi \|_{L^{\infty }(\bbbr )}
%$$
%On the other side, given $\ep'>0$ there exists $x\in [-\12 ,\12 ]$ such that 
%$\| f\|_{L^{\infty }(\bbbr )}\leq |f(x)|+\ep' $ and thus
%$$
%\| fD_{\ep^{-1}}^{1}\varphi \|_{L^{\infty }(\bbbr )}\geq |f(x)\varphi (\ep x)|
%\geq (\| f\|_{L^{\infty }(\bbbr )}-\ep')|\varphi (\ep x)|
%$$
%and so
%$$
%\liminf_{\ep \rightarrow 0} \| fD_{\ep^{-1}}^{1}\varphi \|_{L^{\infty }(\bbbr )}
%\geq (\| f\|_{L^{\infty }(\bbbr )}-\ep')|\varphi (0)|
%=(\| f\|_{L^{\infty }(\bbbr )}-\ep')\| \varphi \|_{L^{\infty }(\bbbr )}
%$$
%for all $\ep'>0$ which finishes the lema.

Now we are ready to proof our main theorem.

\begin{teor} Let $m$ be a bounded continuous function on $\bbbr^{2}$.  Let $0< p_i,q_i<\infty$ for
$i=1,2$,
and $0<p_3,q_3\le\infty$  where
$1/p_1+1/p_2=1/p_3$.

Then $m$ is a multiplier in $(L^{p_1,q_1}(\bbbr)\times L^{p_2,q_2}(\bbbr),L^{p_3,q_3}(\bbbr))$ if and
only if
$(D_{t^{-1}}m)_{t>0}$ restricted to  $\bbbz^{2}$ are
uniformly bounded multipliers in $(L^{p_1,q_1}(\T)\times L^{p_2,q_2}(\T),L^{p_3,q_3}(\T))$,i.e,
denoting $P_t=
P_{(D_{t^{-1}}m)_{k,k'}}
$ where $(D_{t^{-1}}m)_{k,k'}=m(tk,tk')$, there exists $C>0$ such that
$$
\| C_{1}(f,g)\|_{L^{p_{3},q_{3}}(\bbbr )}\leq C\|f\|_{L^{p_{1},q_{1}}(\bbbr )}
\|g\|_{L^{p_{2},q_{2}}(\bbbr )}
$$
for $f,g\in S(\bbbr )$ if and only if there exists $C'>0$ such that
$$
\| P_{t}(f,g)\|_{L^{p_{3},q_{3}}(\toro )}\leq C'\|f\|_{L^{p_{1},q_{1}}(\toro )}
\|g\|_{L^{p_{2},q_{2}}(\toro )}
$$
uniformly in $t>0$ for all trigonometric polynomials $f,g$.
\end{teor}
\Dem
($\Rightarrow $) Let $\varphi =\chi_{[-\12 ,\12 ]}$ and $\psi(x)=\pi^{-1/2}e^{-x^2}$
Let $t>0$ and let
 $f(x)=\sum_{k_{1}\in \bbbz }a_{k_{1}}e^{2\pi ik_{1}x}$ and 
$g(x)=\sum_{k_{2}\in \bbbz }b_{k_{2}}e^{2\pi ik_{2}x}$. 

Since $m$ is continuous we  can write
$$
P_{t}(f,g)(x)=\sum_{k_{1}\in \bbbz }\sum_{k_{2}\in \bbbz }a_{k_{1}}b_{k_{2}}
m(tk_{1},tk_{2})e^{2\pi i(k_{1}+k_{2})x}
$$
$$
=\sum_{k_{1}\in \bbbz }\sum_{k_{2}\in \bbbz }a_{k_{1}}b_{k_{2}}
\lim_{\ep \rightarrow 0}\int_{\bbbr }\int_{\bbbr }
D_{\ep }^{1}\psi (k_{1}-r)D_{\ep }^{1}\psi (k_{2}-s)m(tr,ts)
e^{2\pi i(r+s)x}drds
$$
$$
=\lim_{\ep \rightarrow 0}\int_{\bbbr}\int_{\bbbr}
\sum_{k_{1}\in \bbbz }a_{k_{1}}D_{\ep }^{1}\psi (r-k_{1})\ecu
\sum_{k_{2}\in \bbbz }b_{k_{2}}D_{\ep }^{1}\psi (s-k_{2})\ecu
m(tr,ts)\ecu e^{2\pi i(r+s)x}\ecu drds.
$$
That is \begin{equation}\label{f1}
P_{t}(f,g)(x)
=\lim_{\ep \rightarrow 0}C_{t}(f_{\ep },g_{\ep })(x)
\end{equation}
where 
$$
\hat{f}_{\ep }=\sum_{k_{1}\in \bbbz }a_{k_{1}}T_{k_{1}}D_{\ep }^{1}\psi ,\qquad \hat{g}_{\ep
}=\sum_{k_{2}\in \bbbz
}b_{k_{2}}T_{k_{2}}D_{\ep }^{1}\psi 
$$
or, in other words, 
$$
f_{\ep }(x)=\sum_{k_{1}\in \bbbz }a_{k_{1}}
M_{k_{1}}D_{\ep^{-1}}^{\infty }\check{\psi }(x)
=\sum_{k_{1}\in \bbbz }a_{k_{1}}\check{\psi }(\ep x)e^{2\pi ik_{1}x}
=\check{\psi }(\ep x)f(x),
$$
and similar formula for $g_\ep$. Moreover, this the convergence is 
uniform since 
$$
|P_{t}(f,g)(x)-C_{t}(f_{\ep },g_{\ep })(x)|
$$
$$
\leq \sum_{k_{1}\in \bbbz }\sum_{k_{2}\in \bbbz }|a_{k_{1}}||b_{k_{2}}|
\int_{\bbbr }\int_{\bbbr }|m(tk_{1},tk_{2})-m(t(k_{1}-\ep r),t(k_{2}-\ep s))|
\psi (r)\psi (s)drds
$$
which tends to zero uniformly in $x\in \bbbr $ because the continuity of $m$.

Thus 
\begin{equation}\label{f3}
P_{t}(f,g)= \lim_{n\rightarrow \infty } C_{t}(f_{n},g_{n})
\end{equation}
where $ f_{n}(x)=\check{\psi }(n^{-1}x)f(x) $ and $g_{n}(x)=\check{\psi
}(n^{-1}x)g(x)$ with uniform convergence
and from Lemma \ref{tororealdos} for $k\in \bbbn $ we also have
\begin{equation}\label{f2}
\| P_{t}(f,g)\|_{L^{p_{3},q_3}(\toro )}
=\| P_{t}(f,g)D_{k}^{p_{3}}\varphi \|_{L^{p_{3},q_{3}}(\bbbr )}
\end{equation}

%From Fatou's lemma one has $$\| \liminf_{n\rightarrow \infty }
%\phi_n \|_{L^{p_{3},q_{3}}(\bbbr )}\le  \liminf_{n\rightarrow \infty }\|
%\phi_n \|_{L^{p_{3},q_{3}}(\bbbr )}$$
%for non negative measurable functions.

%Using (\ref{f1}) and (\ref{f2}) together with  the previous observation, we get   

Combining these two facts we write
$$
\| P_{t}(f,g)\|_{L^{p_{3}}(\toro )}
=\| P_{t}(f,g)D_{n}^{p_{3}}\varphi \|_{L^{p_{3}}(\bbbr )}
$$
$$
\leq \| C_{t}(f_{n},g_{n})D_{n}^{p_{3}}\varphi \|_{L^{p_{3}}(\bbbr )}
+\| D_{n^{-1}}( P_{t}(f,g)-C_{t}(f_{n},g_{n}))\varphi \|_{L^{p_{3}}(\bbbr )}
$$

For the first sumand 
\begin{eqnarray*}
\| C_{t}(f_{n},g_{n})
D_{n}^{p_{3}}\varphi \|_{L^{p_{3},q_{3}}(\bbbr )}
&=&\| D^{p_3}_{n}(\varphi D_{n^{-1}}C_{t}(f_{n},g_{n})
) \|_{L^{p_{3},q_{3}}(\bbbr )}
\\
&=&\| \varphi D_{n^{-1}}C_{t}(f_{n},g_{n})
 \|_{L^{p_{3},q_{3}}(\bbbr )}
\\
&\leq& \| D_{n^{-1}}C_{t}(f_{n},g_{n})\|_{L^{p_{3},q_{3}}(\bbbr )}
\| \varphi \|_{L^{\infty }(\bbbr )}\\
&=&n^{-\frac{1}{p_{3}}}
\| C_{t}(f_{n},g_{n})\|_{L^{p_{3},q_{3}}(\bbbr )}
\\
&\leq& n^{-\frac{1}{p_{3}}}
C\|f_{n}\|_{L^{p_{1},q_{1}}(\bbbr )}\|g_{n}\|_{L^{p_{2},q_{2}}(\bbbr )}\\
&=&Cn^{-\frac{1}{p_{1}}}\|f_{n}\|_{L^{p_{1},q_{1}}(\bbbr )}
n^{-\frac{1}{p_{2}}}\|g_{n}\|_{L^{p_{2},q_{2}}(\bbbr )}
\end{eqnarray*}
where, using Lemmas \ref{tororeal1} and \ref{tororeal2}, we know
$$
\lim_{n\rightarrow \infty }n^{-\frac{1}{p_{1}}}\| f_{n}\|_{L^{p_{1},q_{1}}(\bbbr )}
\leq (2^{\frac{r_{1}}{p_{1}}}-1)^{-\frac{1}{r_{1}}}
\| f \|_{L^{p_{1},q_{1}}(\toro )}\|
\check{\psi }\|_{L^{p_{1},r_1}(\bbbr )}
$$
and
$$
\lim_{n\rightarrow \infty }n^{-\frac{1}{p_{2}}}\| g_{n}\|_{L^{p_{2},q_{2}}(\bbbr )}
\leq (2^{\frac{r_{2}}{p_{2}}}-1)^{-\frac{1}{r_{2}}}
\| g \|_{L^{p_{2},q_{2}}(\toro )}\|
\check{\psi }\|_{L^{p_{2},r_2}(\bbbr )}
$$
with $r_{i}=\min(p_{i},q_{i})$ for $i=1,2$.

Thus
$$
\| P_{t}(f,g)\|_{L^{p_{3}}(\toro )}
\leq \lim_{n\rightarrow \infty }\| C_{t}(f_{n},g_{n})
D_{n}^{p_{3}}\varphi \|_{L^{p_{3},q_{3}}(\bbbr )}
$$
$$
+\lim_{n\rightarrow \infty }
\| P_{t}(f,g)-C_{t}(f_{n},g_{n})\|_{L^{\infty }(\bbbr  )} 
=A(p_1,p_2)
\|f\|_{L^{p_{1},q_{1}}(\toro )}\|g\|_{L^{p_{2},q_{2}}(\toro )}
$$
and the proof of this implication is completed.

($\Leftarrow $) Assume $D_{t^{-1}}m$ restricted to  $\bbbz^{2}$ are
uniformly bounded multipliers on $\bbbz^2$ and let $f,g\in S(\bbbr )$ such that $\hat f$ and $\hat g$
have compact support
contained in $K$.

Using  Poisson formula
$$
t\sum_{k_{1}}\hat{f}(tk_{1})e^{2\pi ik_{1}x}
=\sum_{k_{1}}(D_{t}f)\hat{}(k_{1})e^{2\pi ik_{1}x}
=\sum_{k_{1}}D_{t}f(x+k_{1})
=\widetilde{D_{t}f}(x)
$$

Therefore, since $m$ is continuous,  we can write
\begin{eqnarray*}
C_{1}(f,g)(x)&=&\dint_{K\times K}\hat{f}(\xi )\hat{g}(\eta )m(\xi, \eta )
e^{2\pi i(\xi +\eta )x}d\xi d\eta 
\\
&=&\lim_{t\rightarrow 0}t^{2}\sum_{k_{1}}\sum_{k_{2}}\hat{f}(tk_{1})
\hat{g}(tk_{2})m(tk_{1},tk_{2})e^{2\pi it(k_{1}+k_{2})x}\\
&=&\lim_{t\rightarrow 0}P_{t}(\widetilde{D_{t}f},\widetilde{D_{t}g})(tx)
\end{eqnarray*}

Note that
\begin{eqnarray*}
|\{x\in \bbbr: |C_{1}(f,g)(x) |>\lambda\}|&\le &\liminf_{t\to 0} |\{|x|\le t^{-1}/2:
|P_{t}(\widetilde{D_{t}f},\widetilde{D_{t}g})(tx) |>\lambda\}|\\
&\le &\liminf_{t\to 0} t^{-1}|\{|x|\le 1/2:
|P_{t}(\widetilde{D_{t}f},\widetilde{D_{t}g})(x) |>\lambda\}|
\end{eqnarray*}

Therefore, formula  (\ref{equiv}) and Fatou's lemma give
$$||C_1(f,g)||^{p_3}_{L^{p_3,q_3}(\bbbr)}\le C \liminf_{t\to
0}t^{-1}||P_t(\widetilde{D_{t}f},\widetilde{D_{t}g})||^{p_3}_{L^{p_3,q_3}(\T)}.$$

An application of the assumption and Lemma \ref{realtoro} lead to
\begin{eqnarray*}
||C_1(f,g)||_{L^{p_3,q_3}(\bbbr)}&\le&
C \liminf_{t\to 0}
t^{-1/p_3}||\widetilde{D_{t}f}||_{L^{p_1,q_1}(\T)}||\widetilde{D_{t}g}||_{L^{p_2,q_2}(\T)}\\
&\leq &
C||f||_{L^{p_1,q_1}(\bbbr)}||g||_{L^{p_2,q_2}(\bbbr)}.
\end{eqnarray*}

This finishes the proof \qed.

It is known that transference theorems can be extended to symbols more general than continuous (see
\cite{L}, \cite{CW},
\cite{FS}). Actually a bounded measurable function $m_1$ defined on $\bbbr$ is called regulated if
$$
\lim_{\ep\rightarrow 0^+} \frac{1}{2\ep}\int_{-\ep}^{\ep}m_1( x+t)dt=m_1(x)
$$
for all $x\in \bbbr$.

It is pointed out in \cite{L} (see Corollary 2.5 ) that if $m_1$ is regulated and
$\phi$ is non-negative, symmetric, smooth with
compact support and
$\int_{\bbbr}\phi(t)dt=1$ then
$$
\lim_{\ep\rightarrow 0^+} \int_{\bbbr}m_1(x-\ep t)\phi(t)dt=\lim_{\ep\rightarrow 0}
m_1*D^1_\ep\phi(x)=m_1(x)
$$
for all $x\in \bbbr$.

This acually implies that
\begin{equation}\label{reg}
\lim_{\ep\rightarrow 0^+} \int_{\bbbr}m_1(x-\ep t)\psi(t)dt=\lim_{\ep\rightarrow 0} m_1*D^1_\ep
\psi(x)=m_1(x)
\end{equation}
where $\psi$ is non-negative symmetric, smooth and
$\int_{\bbbr}\psi(t)dt=1$.

Indeed, given $\psi$ take non-negative, symmetric, smooth functions $\phi_n$ with
compact support  and
$\int_{\bbbr}\phi_n(t)dt=1$ such that $\lim_{n\to\infty}||\psi-\phi_n||_1=0$ and observe that
$$
 |\int_{\bbbr}(m_1(x-\ep t)-m_1(x))\psi(t)dt|\le
2||m_1||_\infty\int_{\bbbr}|D^1_\ep\psi(t)-D^1_\ep\phi_n(t)|dt$$
$$+
 |\int_{\bbbr}(m_1(x-\ep t)-m_1(x))\phi_n(t)dt|$$
$$=2||m_1||_\infty||\psi-\phi_n||_1
+
 |\int_{\bbbr}(m_1(x-\ep t)-m_1(x))\phi_n(t)dt|.$$

\begin{defi}\label{e1}
Let $G(t,s)=\pi^{-1}e^{-(t^2+s^2)}$. A bounded measurable function $m$ defined on $\bbbr^{2}$ is
$G$-regulated if
$$
\lim_{\ep\to 0} \int_{\bbbr^2}m(x-\ep t,y-\ep s)G(t,s)dtds =\lim_{\ep\to 0} m*D^1_\ep G(x,y)=m(x,y)
$$
for all $(x,y)\in \bbbr^2$.

\end{defi}

A look at the proof of the previous theorem shows that $m$ needs not be continuous but only
$G$-regulated for the argument to
work. 

\begin{teor} Let $m$ be a bounded $G$-regulated function on $\bbbr^{2}$, $0< p_i,q_i<\infty$ for
$i=1,2$ and $0<p_3,q_3\le\infty$  where
$1/p_1+1/p_2=1/p_3$.

If $m$ is a multiplier in $(L^{p_1,q_1}(\bbbr)\times L^{p_2,q_2}(\bbbr),L^{p_3,q_3}(\bbbr))$ then
$m$ restricted to  $\bbbz^{2}$ is
a bounded multiplier in $(L^{p_1,q_1}(\T)\times L^{p_2,q_2}(\T),L^{p_3,q_3}(\T))$.
\end{teor}
Now we can apply this result to transfer results for the bilinear Hilbert transform because of the
following remark.

\begin{nota} If $m_1$ be a regulated function defined in $\R$ then $m_\alpha(x,y)=m_1(x+\alpha y)$ is
$G$-regulated in
$\bbbr^2$.

In particular, $m(x,y)=sign(x+\alpha y)$ is $G$-regulated.
\end{nota}
Indeed, observe that
\begin{eqnarray*}
\int_{\bbbr^2}m_1(x- t+\alpha (y- s))D^1_\ep G(t,s)dtds&=&\int_{\bbbr}\int_{\bbbr} m_1(x+\alpha y-
\ep(t+\alpha
s))G(t,s)dtds\\
&=&\int_{\bbbr}m_1(x+\alpha y- \ep t
)\big(\int_{\bbbr}G(t-\alpha s,s)ds\big)dt\\
&=&\int_{\bbbr}m_1(x+\alpha y- \ep t
)\psi_{\al }(t)dt
\end{eqnarray*}
where $\psi_{\al }(t)=\int_{\bbbr}G(t-\alpha s,s)ds$. 
Hence we have, from (\ref{reg}), that  
$$\lim_{\ep\to 0}\int_{\bbbr^2}m_\alpha(x-t,y-s)D^1_\ep G(t,s)dtds=m_\alpha(x,y).$$

\end{document}